\documentclass[11pt]{article}

\usepackage{pgfpages}
\usepackage{amsmath,amssymb,verbatim}
\usepackage{listings}
\usepackage{epsfig}
\usepackage[latin1]{inputenc}

\newtheorem{thm}{Theorem}
\newtheorem{defn}{Definition}
\newtheorem{rem}{Remark}

\def\blacksquare{\hbox{\vrule width 7pt height 7pt depth 0pt}}

\begin{document}

\tolerance=5000
\setcounter{page}{1}

\centerline{{\Large{\bf Initial-boundary-value problems for the one-}}}
\vspace{0.1cm}

\centerline{{\Large{\bf dimensional time-fractional diffusion equation}}}

\vspace{0.3cm}
\centerline{{\bf Yuri Luchko}}
\vspace{0.3cm}

\centerline{{Department of Mathematics}}

\centerline{{Beuth Technical University of Applied Sciences Berlin}}

\centerline{{Luxemburger Str. 10, 13353 Berlin,\ Germany}}

\vspace{0.2cm}

\begin{abstract}
\noindent
In this paper, some initial-boundary-value problems for the time-fractional diffusion equation 
are first considered in open bounded n-dimensional domains. In particular, 
the maximum principle well-known for the PDEs of elliptic and parabolic types is extended for the time-fractional diffusion equation.  In its turn, the maximum principle is used to show 
 uniqueness of solution to the   
initial-boundary-value problems for the time-fractional  diffusion equation.  The generalized solution in sense of Vladimirov is then constructed in form of a Fourier series with respect to the eigenfunctions of a certain Sturm-Liouville eigenvalue problem.  For the one-dimensional time-fractional
diffusion equation  
$$
(D_t^{\alpha} u)(t) = \frac{\partial}{\partial x}\left( p(x)
\frac{\partial u}{\partial x}\right) -q(x)\, u + F(x,t),\ \ 
x\in (0,l),\ t\in (0,T)
$$
the generalized solution to the   
initial-boundary-value problem with the Dirichlet boundary conditions is shown to be a solution in the classical sense. Properties of the solution  are investigated including its smoothness and asymptotics for some special cases of the source function. 
\end{abstract}

\vspace{0.2cm}



\section{Introduction}

During the last few decades, PDEs of fractional order  have been successfully used for
modeling of many  processes and systems that stimulated development of their mathematical theory. This
process is still going on, but one can already recognize some new  insights and correlations
between different branches of mathematics and natural sciences that have
emerged within the framework
of FC models. The literature related both to theory of PDEs of fractional order and to their applications is very 	 extensive. For the theory of fractional differential equations we refer e.g. to the books \cite{Dit10}, \cite{Kil06}, and \cite{Pod99}. As to the applications, we mention the book \cite{Uch08} that is completely devoted to different applications of  FC  in physics, chemistry, technique, astrophysics, etc. and contains several dozens of interesting case studies. The recent book \cite{Her11} deals with applications of fractional differential equations in classical mechanics, group theory, quantum
mechanics, nuclear physics, hadron spectroscopy, and quantum field theory. For other interesting FC models we refer the reader to \cite{Dub}, \cite{FDL}, \cite{Hil00}, \cite{Kla08}, \cite{LucPun11}, \cite{Mag04}, \cite{Mai10}, \cite{Mai96}, \cite{Met04} to mention only few of many recent publications. 

In applications, relevant processes like e.g. anomalous diffusion or wave propagation in complex systems mainly run  within some bounded domains in space that corresponds to initial-boundary-value problems for the PDEs of fractional order, which model the processes under consideration. Both analytical and numerical investigation of these problems  thus become an important task in Fractional Calculus. 

In the author's papers \cite{Luc11}-\cite{Luc09_2}, the  time-fractional diffusion equation 
$$
(D_t^{\alpha} u)(t) = \mbox{div}(p(x)\,\mbox{grad}\,u) - q(x)u\, + F(x,t), \ \ 0<\alpha\le 1 
$$
with
the Caputo fractional derivative $D_t^{\alpha}$ 
with respect to the time variable $t$  
  along with some of its important generalizations (multi-term equation and equation of distributed order) and  some initial and boundary conditions has been considered in open bounded n-dimensional domains. 
In particular, 
the maximum principle well-known for the PDEs of elliptic and parabolic types was extended for the time-fractional diffusion equation as well as for the multi-term time-fractional diffusion equation and for the time-fractional diffusion equation of distributed order. The validity  of the maximum principle is based on an important extremum principle for the Caputo fractional derivative. In its turn, the maximum principle was used to show 
 uniqueness of solution to some   
initial-boundary-value problems for the time-fractional  diffusion equations.  This solution - if it exists - continuously depends on the problem data, i.e., on the source function and on the initial and boundary conditions.  

As to the existence of solution, in the author's papers \cite{Luc11}, \cite{Luc10_1}, \cite{Luc09_2}, the Fourier method of variables separation was used to obtain a formal solution in form of a Fourier series with respect to the eigenfunctions of a certain Sturm-Liouville eigenvalue problem. Under certain conditions, the formal solution was shown to be a generalized solution in sense of Vladimirov, i.e., a continuous function that satisfies the PDE of fractional order in generalized sense (see \cite{Vla}). To prove that this generalized solution is a solution in the classical sense, i.e., it is at least twice differentiable function with respect to the spatial variables and $\alpha$-differentiable with respect to the time variable turned out to be a difficult undertaking that is still not resolved in general case.

In this paper, this problem is solved for the initial-boundary-value problems with the Dirichlet boundary conditions
for the one-dimensional time-frac\-ti\-on\-al
diffusion equation. Under certain conditions, the unique solution to this problem can be represented via a Fourier series with respect to the eigenfunctions of a corresponding Sturm-Liouville eigenvalue problem.

For other methods and results in the theory of initial-boundary-value problems for the PDEs of fractional order we refer e.g. to  \cite{Baz98}, \cite{Gor00}, \cite{Kem11}, \cite{Mee09}, \cite{Met00a}, \cite{Psk05} to mention only few of many recent publications. 

The rest of the paper is organized as follows. In the 2nd section,  definitions and notions used in the paper  along with the problem formulation are given. In the 3rd section, main results for the initial-boundary-value problems for the time-fractional diffusion equation 
in open bounded n-dimensional domains are short presented. In particular, 
the maximum principle well-known for the PDEs of elliptic and parabolic types is extended for the time-fractional diffusion equation and used to show 
 uniqueness of solution to the   
initial-boundary-value problems for the time-fractional  diffusion equation.  The generalized solution in sense of Vladimirov is then constructed in form of a Fourier series with respect to the eigenfunctions of a certain Sturm-Liouville eigenvalue problem. The 4rd section deals with the one-dimensional time-fractional
diffusion equation.   In this case, the generalized solution is shown to be a solution in the classical sense. Properties of the solution  are investigated including its smoothness and asymptotics for some special cases of the source function $F$. Finally, the last section contains some conclusions and open problems. 

\section{Problem formulation}

In this paper, we deal with the time-fractional diffusion equation
\begin{equation}
\label{eq}
(D_t^{\alpha} u)(t) = - L(u) + F(x,t),
\end{equation}
$$
0<\alpha\le 1,
(x,t)\in \Omega_T:= G \times (0,T),\ G \subset  R^n,
$$
and especially with its one-dimensional case ($n=1$),
where
$$
L(u):= - \mbox{div}(p(x)\,\mbox{grad}\,u) + q(x)u,
$$
\begin{equation}
\label{pq}
p\in C^1(\bar{G}), q\in C(\bar{G}), \ p(x)>0, q(x)\ge 0, x \in \bar{G},
\end{equation}
the fractional derivative $D_t^{\alpha}$ is defined in the Caputo sense
\begin{equation}
\label{fd}
(D^{\alpha} f)(t):= (I^{1-\alpha} f^\prime )(t), \ 0<\alpha\le 1,
\end{equation}
$I^\alpha$ being the Riemann-Liouville fractional integral
$$
(I^{\alpha} f)(t):= 
\begin{cases}
{1\over \Gamma(\alpha)} \int_0^{t} (t-\tau)^{\alpha -1} f(\tau)\, d\tau, \  0<\alpha <1, \\
f(t),\ \alpha =0,
\end{cases}
$$
and the domain $G$ with the boundary $S$ is open and bounded in
$R^n$. For the theory of the fractional integrals and derivatives
the reader is referred  e.g. to \cite{Kil06} or \cite{Pod99}.

For $\alpha=1$, the equation (\ref{eq}) is reduced to a linear second-order
parabolic PDE. The theory of this equation is well-known, so that the main focus
in the paper is on the case $0<\alpha <1$.
For the equation (\ref{eq}), the initial-boundary-value problem
\begin{equation}
\label{ic}
u\big|_{t=0}=u_0(x),\ x\in \bar{G},
\end{equation}
\begin{equation}
\label{bc}
u\big|_{S}=v(x,t),\ (x,t)\in S\times [0,T]
\end{equation}
is considered. A solution of the problem (\ref{eq}),
(\ref{ic}), (\ref{bc}) is called a function $u=u(x,t)$ defined in
the domain $\bar{\Omega}_T:= \bar{G} \times [0,T]$ that belongs to
the space $C(\bar{\Omega}_T) \cap W^1_t((0,T])\cap C^2_x(G)$ and
satisfies both the equation (\ref{eq}) and the initial and
boundary conditions (\ref{ic})-(\ref{bc}). By $W^1_t((0,T])$,  the
space of the functions $f \in C^1((0,T])$ such that $f^\prime \in
L((0,T))$ is denoted. If the problem (\ref{eq}), (\ref{ic}),
(\ref{bc}) possesses a  solution, then the functions $F$,
$u_0$ and $v$ given in the problem have to belong to the spaces
$C(\Omega_T),\ C(\bar{G})$ and $C(S\times [0,T])$, respectively.
In the further discussions,  these inclusions are always supposed
to be valid.

\section{n-dimensional time-fractional diffusion equation}

In this section, a short survey of results for initial-boundary-value problems (\ref{ic})-(\ref{bc}) for the n-dimensional time-fractional diffusion equation (\ref{eq}) is presented. For the proofs of results we refer the readers to \cite{Luc10_1} and \cite{Luc09_1}. 

We start with  uniqueness
of solution of the problem (\ref{eq}), (\ref{ic}), (\ref{bc}). The
main component of the uniqueness proof is an
appropriate maximum principle for the equation (\ref{eq}).
%
%
%

\begin{thm}
\label{t2}
Let a function $u \in C(\bar{\Omega}_T) \cap W^1_t((0,T])\cap C^2_x(G)$ be
 a   solution of the
generalized time-fractional diffusion equation (\ref{eq}) in the domain $\Omega_T$
and $F(x,t)\le 0,\ (x,t)\in \Omega_T$.

Then either
$u(x,t)\le 0,\ (x,t)\in \bar{\Omega}_T$ or
the function $u$
attains its positive maximum on the part $S_G^T:=
(\bar{G}\times\{0\}) \cup  (S\times[0,T])$
of the boundary
of the domain
$\Omega_T$, i.e.,
\begin{equation}
\label{max}
u(x,t)\le \max \{0, \max\limits_{ (x,t)\in S_G^T } u(x,t)\},\ \forall (x,t)\in \bar{\Omega}_T.
\end{equation}
\end{thm}

Similarly to the case of the PDEs of  parabolic type ($\alpha =1$),
an appropriate minimum principle is valid, too.

%

The maximum and minimum principles can be applied to show that the
 problem (\ref{eq}), (\ref{ic})-(\ref{bc}) possesses at most one  solution
and this solution - if it exists - continuously depends on
the data given in the problem.

\begin{thm}
\label{t5}
The initial-boundary-value-problem (\ref{ic})-(\ref{bc})  for the time-fractional diffusion equation (\ref{eq})
 possesses at most one  solution.
This solution continuously depends on the data given in the problem
in the sense that if
$$
\|F-\tilde{F}\|_{C(\bar{\Omega}_T)}\le \epsilon,
\| u_0-\tilde{u}_0\|_{C(\bar{G})}\le \epsilon_0,\
\| v-\tilde{v}\|_{C(S\times [0,T])} \le \epsilon_1,
$$
and $u$ and $\tilde{u}$ are the  solutions of the
problem (\ref{eq}), (\ref{ic})-(\ref{bc}) with the source functions
$F$ and $\tilde{F}$, the initial conditions $u_0$ and $\tilde{u_0}$,
and the boundary conditions $v$ and $\tilde{v}$, respectively, then the norm estimate
\begin{equation}
\label{est2}
\| u-\tilde{u}\|_{C(\bar{\Omega}_T)}
  \le \max\{ \epsilon_0, \epsilon_1\}\, + {T^\alpha
\over \Gamma(1+\alpha)}\, \epsilon
\end{equation}
for the solutions $u$ and $\tilde{u}$ holds true.
\end{thm}

To tackle the
problem of existence of the solution, the notion of the
generalized solution is introduced (see \cite{Vla} for the case $\alpha =1$).

\begin{defn} 
\label{d1}
Let $F_k  \in C(\bar{\Omega}_T),\ u_{0k}\in C(\bar{G})$ and
$v_k\in C(S\times [0,T]),\ k=1,2,\dots$ be  the sequences of functions that
satisfy the following conditions:

\noindent
1) there exist the functions $F$, $u_0$, and $v$, such that
\begin{equation}
\label{Fk}
\|F_k- F \|_{C(\bar{\Omega}_T)} \to 0 \ \ \mbox{as} \ \ k\to +\infty,
\end{equation}
\begin{equation}
\label{u0k}
\| u_{0k}- u_0\|_{C(\bar{G})} \to 0 \ \ \mbox{as} \ \ k\to +\infty, \
\end{equation}
\begin{equation}
\label{vk}
\| v_k- {v}\|_{C(S\times [0,T])} \to 0 \ \ \mbox{as} \ \ k\to +\infty,
\end{equation}

\noindent
2) for any $k=1,2,\dots$ there exists a solution $u_k$
 of the initial-boundary-value
problem
\begin{equation}
\label{ick}
u_k\big|_{t=0}=u_{0k}(x),\ x\in \bar{G},
\end{equation}
\begin{equation}
\label{bck}
u_k\big|_{S}=v_k(x,t),\ (x,t)\in S\times [0,T],
\end{equation}
for the time-fractional diffusion equation
\begin{equation}
\label{eqk}
(D_t^{\alpha} u_k)(t) = -L(u_k) + F_k(x,t).
\end{equation}
Suppose, there exists a function $u\in C(\bar{\Omega}_T)$ such that
\begin{equation}
\label{uk}
\| u_{k}- u\|_{C(\bar{G})} \to 0 \ \ \mbox{as} \ \ k\to +\infty.
\end{equation}
The function $u$ is called a generalized solution of the problem (\ref{eq}), (\ref{ic})-(\ref{bc}).
\end{defn}

A generalized solution of the problem (\ref{eq}),
(\ref{ic})-(\ref{bc}) is a continuous function, not a generalized
one. Still, the generalized solution is not required to be from
the functional space $C(\bar{\Omega}_T) \cap W^1_t((0,T])\cap
C^2_x(G)$, where the  solution has to belong to.

It follows from Definition 1 that if the problem
(\ref{eq}), (\ref{ic})-(\ref{bc}) possesses a  solution
then this solution is a generalized solution of the problem, too.
In this sense, Definition 1 extends the notion of 
 solution of the problem (\ref{eq}),
(\ref{ic})-(\ref{bc}). This extension is needed to get some
existence results. But of course one does not want to lose 
uniqueness of the solution. The following result is valid:

\begin{thm}
\label{t6}
The problem (\ref{eq}), (\ref{ic})-(\ref{bc}) possesses at most one generalized solution
in the sense of Definition 1.
The generalized solution - if it exists - continuously depends on the data given in the problem
in the sense of the estimate (\ref{est2}).
\end{thm}

Contrary to the situation with  solution of the
problem (\ref{eq}), (\ref{ic})-(\ref{bc}), existence of the
generalized solution
 can be shown in the general case under some standard
restrictions on the problem data and the boundary $S$ of the domain $G$.
In this section,  existence of the generalized solution of the problem
\begin{equation}
\label{eq0}
(D_t^{\alpha} u)(t) = -L(u),
\end{equation}
\begin{equation}
\label{ic0}
u\big|_{t=0}=u_{0}(x),\ x\in \bar{G},
\end{equation}
\begin{equation}
\label{bc0}
u\big|_{S}=0,\ (x,t)\in S\times [0,T]
\end{equation}
is considered to demonstrate the technique that can be used with
the appropriate standard modifications in the general case, too (see the case $n=1$ in the next section).
The generalized solution of the problem (\ref{eq0})-(\ref{bc0})
can be constructed in an analytical form by using the Fourier
method of the variables separation. Let us look for a particular
solution $u$ of the equation (\ref{eq0}) in the form
\begin{equation}
\label{tx}
u(x,t) = T(t)\, X(x),\ \ (x,t)\in \bar{\Omega}_T,
\end{equation}
that satisfies the boundary condition (\ref{bc0}). Substitution of the
function (\ref{tx}) into the equation (\ref{eq0}) and separation
of the variables lead to the equation
\begin{equation}
\label{tx2}
{(D_t^{\alpha} T)(t) \over T(t)} = -{L(X) \over X(x)}\, = \, -\lambda,
\end{equation}
$\lambda$ being a constant not depending on the variables $t$ and
$x$. The last equation, together with the boundary condition
(\ref{bc0}), is equivalent to the fractional differential equation
\begin{equation}
\label{tx3}
(D_t^{\alpha} T)(t) + \lambda T(t)  = 0
\end{equation}
and the eigenvalue problem
\begin{equation}
\label{tx4}
L(X) = \lambda\, X,
\end{equation}
\begin{equation}
\label{tx5}
X\big|_{S}=0,\ x\in S
\end{equation}
for the operator $L$. Due to the condition (\ref{pq}), the
operator $L$ is a positive definite and self-adjoint linear
operator. The theory of the eigenvalue problems 
for such operators is well-known (see e.g.
\cite{Vla}). In particular, the eigenvalue problem  (\ref{tx4}) -
(\ref{tx5}) has a
counted number of the positive eigenvalues $0<\lambda_1\le
\lambda_2\le \dots$ with the finite multiplicity and - if the
boundary $S$ of $G$ is a smooth surface -  any function $f\in
{\cal M}_L$ can be represented through its Fourier series in the
form
\begin{equation}
\label{f1}
f(x) = \sum_{i=1}^\infty (f, X_i)\, X_i(x),
\end{equation}
where $X_i\in {\cal M}_L$ are the eigenfunctions corresponding to the eigenvalues $\lambda_i$:
\begin{equation}
\label{f2}
L(X_i) = \lambda_i\, X_i,\ i=1,2,\dots .
\end{equation}
By  ${\cal M}_L$, the space of the functions $f$
 that satisfy the boundary condition (\ref{tx5}) and
the inclusions $f\in C^1(\bar{\Omega}_T)\cap C^2(G)$, $L(f)\in L^2(G)$ is denoted.

The solution of the fractional differential equation (\ref{tx3}) with $\lambda = \lambda_i,\ i=1,
2,\dots$ has the form (see e.g.  \cite{Luc99A},  \cite{Luc99B})
\begin{equation}
\label{f3}
T_i(t) = c_i E_{\alpha}(-\lambda_i t^\alpha),
\end{equation}
$E_\alpha$ being the Mittag-Leffler function defined by
\begin{equation}
\label{f4}
E_{\alpha}(z):= \sum_{k=1}^\infty {z^k\over \Gamma(\alpha\,k +1)}.
\end{equation}
Any of the functions
\begin{equation}
\label{f5-1}
u_i(x,t) = c_i E_{\alpha}(-\lambda_i t^\alpha)\, X_i(x),\ i=1,2,\dots
\end{equation}
and thus the finite sums
\begin{equation}
\label{f5}
u_k(x,t) = \sum_{i=1}^k c_i E_{\alpha}(-\lambda_i t^\alpha)\, X_i(x),\ k=1,2\dots
\end{equation}
satisfy both the equation (\ref{eq0}) and the boundary condition (\ref{bc0}).
To construct a function that satisfies the initial condition (\ref{ic0}), 
the notion of a formal solution is introduced.

\begin{defn}
A formal solution of the problem (\ref{eq0})-(\ref{bc0}) is called the Fourier
series in the form
\begin{equation}
\label{f6}
u(x,t) = \sum_{i=1}^\infty (u_0, X_i) E_{\alpha}(-\lambda_i t^\alpha)\, X_i(x),
\end{equation}
 $X_i,\ i=1,2,...$ being the eigenfunctions corresponding to the eigenvalues $\lambda_i$
of the eigenvalue problem (\ref{tx4}) - (\ref{tx5}) and $(f,g)$ the standard scalar product in $L^2(0,l)$. 
\end{defn}

Under certain conditions, the formal solution (\ref{f6}) can be proved to be
 the generalized solution of the problem (\ref{eq0})-(\ref{bc0}).

\begin{thm}
\label{t7} 
Let the function $u_0$ in the initial condition
(\ref{ic0}) be from the space ${\cal M}_L$. Then the formal
solution (\ref{f6}) of the problem (\ref{eq0})-(\ref{bc0}) is its
generalized solution.
\end{thm}

To prove that this generalized solution is a solution in the classical sense, i.e., it is at least twice differentiable function with respect to the spatial variables and $\alpha$-differentiable with respect to the time variable turned out to be a difficult undertaking that is still not resolved in general case.

In the next section, this problem is solved for the initial-boundary-value problem for the one-dimensional time-fractional
diffusion equation.

\begin{rem}
The method presented in this section  can be applied  with some small modifications to the case of the problem  (\ref{eq0})-(\ref{bc0}) 
on an infinite domain $\Omega = G \times (0,\infty),\ G \subset  R^n$, too.
\end{rem}

\section{One-dimensional time-fractional diffusion equation}

All results presented in the previous section are valid for the n-dimensional time-fractional diffusion equation and thus for the one-dimensional case, too. The main result of this section is that  the generalized solution of the initial-boundary-value problems for the one-dimensional  time-fractional diffusion equation we received in the previous section can be interpreted as solution in the classical sense under some suitable conditions. 

We start with the initial-boundary-value problem for the inhomogeneous one-dimensional  time-fractional diffusion equation 
\begin{equation}
\label{4-eq}
(D_t^{\alpha} u)(t) = {\partial \over \partial x}\left( p(x)
{\partial u\over \partial x}\right) -q(x)\, u + F(x,t),\ \ 
x\in (0,l),\ t\in (0,T)
\end{equation}
with the inhomogeneous Dirichlet boundary conditions
\begin{equation}
\label{4-bc}
u(0,t)=\phi_1(t),\ u(l,t)=\phi_2(t),\ \ 0\le t \le T,
\end{equation}
along with the initial condition
\begin{equation}
\label{4-ic}
u(x,0)=u_0(x),\ \ 0\le x \le l,
\end{equation}
and the conjugate condition
\begin{equation}
\label{4-cc}
u_0(0)=\phi_1(0).
\end{equation}
In the following, we always suppose that the conjugate condition (\ref{4-cc}) holds true. Like in the case of the initial-boundary-value problems for the PDEs of the parabolic type ($\alpha =1$ in (\ref{4-eq})), the problem (\ref{4-eq})-(\ref{4-cc}) can be transformed to a problem with the homogeneous boundary conditions. Indeed, let us introduce an auxiliary function $v$ according to the formula
\begin{equation}
\label{4-v}
v(x,t)=u(x,t) +\frac{x}{l}(\phi_1(t)-\phi_2(t))-\phi_1(t). 
\end{equation}
Then the inhomogeneous boundary conditions (\ref{4-bc}) for the unknown function $u$ are transformed to the homogeneous boundary conditions for the function $v$:
$$
v(0,t)=0,\ v(l,t)=0,\ \ 0\le t \le T.
$$
The form of the time-fractional diffusion equation (\ref{4-eq}) for the new unknown function $v$ 
$$
(D_t^{\alpha} v)(t) = {\partial \over \partial x}\left( p(x)
{\partial v\over \partial x}\right) -q(x)\, v + F_1(x,t),\ \ 
x\in (0,l),\ t\in (0,T)
$$
with
$$
F_1(x,t) = F(x,t) +\frac{x}{l}(D_t^{\alpha} (\phi_1-\phi_2))(t) - (D_t^{\alpha} \phi_1)(t) - \frac{\phi_1(t)-\phi_2(t)}{l}p^\prime(x) +
$$
$$
q(x)\left(\frac{x}{l}(\phi_1(t)-\phi_2(t))-\phi_1(t)\right)
$$
as well as the form of the initial condition (\ref{4-ic})
$$
v(x,0) = u_0(x) + \frac{x}{l}(\phi_1(0)-\phi_2(0))-\phi_1(0)
$$
remain the same as for the original problem (\ref{4-eq})-(\ref{4-ic}). 

Without loss of generality, we can thus consider the initial-boundary-value problem (\ref{4-eq})-(\ref{4-ic}) with the homogeneous boundary conditions, i.e. with
$$
\phi_1(t)\equiv \phi_2(t)\equiv 0,\ \ 0\le t \le T.
$$
Let us start with the homogeneous equation (\ref{4-eq}), i.e. we first suppose that $F(x,t)\equiv 0,\ 0<x<l,\ 0<t<T$. In the previous section, it was shown (see Theorem 4) that if $u_0 \in {\cal M}_L$ then the formal solution (\ref{f6}) is the generalized solution to the problem (\ref{4-eq})-(\ref{4-ic}). We prove now that (\ref{f6}) is a solution in the classical sense under some additional conditions. 

\begin{thm}
Let $u_0 \in {\cal M}_L$ and $L(u_0) \in {\cal M}_L$.  Then solution of
the initial-boundary-value problem
\begin{equation}
\label{4-ic-h}
u\big|_{t=0}=u_{0}(x),\ 0\le x \le l,
\end{equation}
\begin{equation}
\label{4-bc-h}
u(0,t)=u(l,t)=0,\ 0\le t\le T
\end{equation}
for the one-dimensional time-fractional diffusion equation
\begin{equation}
\label{4-eq-h}
(D_t^{\alpha} u)(t) = {\partial \over \partial x}\left( p(x)
{\partial u\over \partial x}\right) -q(x)\, u
\end{equation}
exists, is unique, and is given by the formula 
\begin{equation}
\label{4-f6}
u(x,t) = \sum_{i=1}^\infty (u_0, X_i) E_{\alpha}(-\lambda_i t^\alpha)\, X_i(x),
\end{equation}
 $X_i,\ i=1,2,...$ being the eigenfunctions corresponding to the eigenvalues $\lambda_i$
of the eigenvalue problem (\ref{tx4}) - (\ref{tx5}).
\end{thm}

Proof: The uniqueness of solution follows from Theorem 2. Theorem 4 states that (\ref{4-f6}) is a generalized solution to the problem (\ref{4-ic-h})-(\ref{4-eq-h}). We show now that (\ref{4-f6}) is a solution in the classical sense, i.e. it is at least twice differentiable function with respect to the spatial variable and $\alpha$-differentiable with respect to the time variable. Indeed, let us differentiate   the series from the left-hand side of the formula (\ref{4-f6}) with respect to the spatial variable term by term and construct a series of derivatives:
\begin{equation}
\label{4-f6-1}
\sum_{i=1}^\infty c_i E_{\alpha}(-\lambda_i t^\alpha)\, X_i^\prime(x),\ \ c_i = (u_0, X_i). 
\end{equation}
Because of the inclusion $u_0\in {\cal M}_L$, the function $L(u_0)$ is from the space $L^2(0,l)$. We can thus determine the Fourier coefficients of $L(u_0)$ with respect to the eigenfunctions $X_i,\ i=1,2,\dots$:
\begin{equation}
\label{4-f6-2}
(L(u_0),X_i)= (u_0, L(X_i)) = \lambda_i (u_0, X_i) = \lambda_ic_i. 
\end{equation}
In the formula (\ref{4-f6-2}), we used the fact that $L$ is a positive definite and self-adjoint operator. The Parseval equality for $L(u_0)$ leads to the relation
\begin{equation}
\label{4-f6-3}
\sum_{i=1}^\infty \lambda_i^2c_i^2 = \left\|L(u_0)\right\|^2<+\infty. 
\end{equation}
It is well-known that the system $X_i,\ i=1,2,\dots$ of the eigenfunctions 
of the eigenvalue problem (\ref{tx4}) - (\ref{tx5}) is a complete orthonormal system in $L^2(0,l)$ and the series 
\begin{equation}
\label{4-f6-4}
\sum_{i=1}^\infty \frac{|X_i(x)|^2}{\lambda_i},\ \ \ 
\sum_{i=1}^\infty \frac{|X_i^\prime(x)|^2}{\lambda_i^2},\ \ \ 
\sum_{i=1}^\infty \frac{|X_i^{\prime\prime}(x)|^2}{\lambda_i^3}
\end{equation}
are uniformly convergent on the interval $[0,l]$ (see e.g. \cite{Vla}).

Using this fact, the estimate (\ref{4-f6-3}), and applying the Cauchy-Bunjakovski inequality we can now prove the uniform convergence of the series (\ref{4-f6-1}):
$$
\sum_{i=1}^\infty |c_i E_{\alpha}(-\lambda_i t^\alpha)\, X_i^\prime(x)| \le 
M \sum_{i=1}^\infty |c_i\, \lambda_i|\left| \frac{X_i^\prime(x)}{\lambda_i}\right| \le 
$$
$$
M \left( \sum_{i=1}^\infty c_i^2\, \lambda_i^2\right)^{\frac{1}{2}} \left( \sum_{i=1}^\infty \frac{|X_i^\prime(x)|^2}{\lambda_i^2}\right)^{\frac{1}{2}} 
\le M \left\|L(u_0)\right\| \left( \sum_{i=1}^\infty \frac{|X_i^\prime(x)|^2}{\lambda_i^2}\right)^{\frac{1}{2}}.
$$
In the first inequality, we used 
 the estimate (see e.g. \cite{Pod99})
\begin{equation}
\label{MLe}
|E_\alpha (-x)|\le {M\over 1+x}\le M,\ 0\le x,\ 0<\alpha<1.
\end{equation}
Because the series (\ref{4-f6-1}) is uniformly convergent one, too (see \cite{Luc10_1}), we thus have proved that the generalized solution (\ref{4-f6}) is a $C^1_x(0,l)$ function and the relation 
$$
\frac{\partial u}{\partial x} = \sum_{i=1}^\infty c_i E_{\alpha}(-\lambda_i t^\alpha)\, X_i^\prime(x)
$$
holds true. The same arguments can be applied to show that (\ref{4-f6}) is a $C^2_x(0,l)$ function and the relation 
\begin{equation}
\label{4-f6-5}
\frac{\partial^2 u}{\partial x^2} = \sum_{i=1}^\infty c_i E_{\alpha}(-\lambda_i t^\alpha)\, X_i^{\prime\prime}(x)
\end{equation}
holds true. Indeed, due to the inclusion $L(u_0)\in {\cal M}_L$ we get $L^2(u_0)\in L^2(0,l)$ and
$$
(L^2(u_0),X_i)= (u_0, L^2(X_i)) = \lambda_i^2 (u_0, X_i) = \lambda_i^2c_i,\ \ i=1,2,\dots .
$$
Then 
$$
\sum_{i=1}^\infty \lambda_i^4c_i^2 = \left\|L^2(u_0)\right\|^2<+\infty 
$$
and we can use the uniform convergence of the 3rd of the series in (\ref{4-f6-4}) and again apply the Cauchy-Bunjakovski inequality to prove the formula (\ref{4-f6-5}):
$$
\sum_{i=1}^\infty |c_i E_{\alpha}(-\lambda_i t^\alpha)\, X_i^{\prime\prime}(x)| \le 
M \sum_{i=1}^\infty |c_i\, \lambda_i^{\frac{3}{2}}|\left| \frac{X_i^{\prime\prime}(x)}{\lambda_i^{\frac{3}{2}}}\right| \le 
$$
$$
M_1 \left( \sum_{i=1}^\infty c_i^2\, \lambda_i^4\right)^{\frac{1}{2}} \left( \sum_{i=1}^\infty \frac{|X_i^{\prime\prime}(x)|^2}{\lambda_i^3}\right)^{\frac{1}{2}} 
\le M_1 \left\|L^2(u_0)\right\| \left( \sum_{i=1}^\infty \frac{|X_i^{\prime\prime}(x)|^2}{\lambda_i^3}\right)^{\frac{1}{2}}.
$$
Note, that in the 2nd inequality we used the inequality $\lambda_i^3 < C\lambda_i^4,\ i=1,2,\dots$ with a certain constant $C$ that holds true because of $\lim_{i\to \infty} \lambda_i = +\infty$.  

To prove that the series given by (\ref{4-f6}) is $\alpha$-differentiable with respect to the time variable we use the following result from \cite{Samko}: 

\noindent
Let for a sequence of functions $f_i,\ i=1,2,\dots$ defined on an interval $(a,b]$ following conditions be fulfilled:

\noindent
1) For a given $\alpha >0$ there exists fractional derivatives $(D^\alpha f_i)(t),\ t\in (a,b]$ for all $i=1,2,\dots$,

\noindent
2) both the series $\sum_{i=1}^\infty f_i(t)$ and the series $\sum_{i=1}^\infty (D^\alpha f_i)(t)$ are uniformly convergent on the interval $[a+\epsilon,b]$ for any $\epsilon >0$.

\noindent
Then the function defined by the series $\sum_{i=1}^\infty f_i(t)$ is $\alpha$-differentiable and the relation
$$
(D^\alpha \sum_{i=1}^\infty f_i)(t) = \sum_{i=1}^\infty (D^\alpha f_i)(t),\ \ a < t < b
$$
holds true. 

Now let us apply the fractional derivative of order $\alpha$ with respect to the time variable to the series from the left-hand side of the formula (\ref{4-f6})  term by term and construct a series of derivatives:
\begin{equation}
\label{4-f6-6}
\sum_{i=1}^\infty c_i (D^\alpha E_{\alpha}(-\lambda_i t^\alpha))(t)\, X_i(x)\ = \ - \sum_{i=1}^\infty c_i \lambda_i E_{\alpha}(-\lambda_i t^\alpha)\, X_i(x).
\end{equation}

The formula (\ref{4-f6-6}) follows from the fact that the Mittag-Leffler function $E_{\alpha}(-\lambda_i t^\alpha)$ is a solution of the fractional differential equation (\ref{tx3}) with $\lambda = \lambda_i,\ i=1,
2,\dots$ or can be directly verified using the well-known formula ($0< \alpha \le 1,\ 0<\beta$) 
\begin{equation}
\label{d-c}
(D^{\alpha} \tau^\beta )(t) = {\Gamma(1+\beta)\over \Gamma(1-\alpha +\beta)} t^{\beta-\alpha}.
\end{equation}
Then we apply the Cauchy-Bunjakovski inequality to the series 
(\ref{4-f6-6}) and get the estimates
$$
\sum_{i=1}^\infty |-c_i \lambda_i E_{\alpha}(-\lambda_i t^\alpha)\, X_i(x)| \le 
M \sum_{i=1}^\infty |c_i\, \lambda_i^{\frac{3}{2}}|\left| \frac{X_i(x)}{\lambda_i^{\frac{1}{2}}}\right| \le 
$$
$$
M_1 \left( \sum_{i=1}^\infty c_i^2\, \lambda_i^4\right)^{\frac{1}{2}} \left( \sum_{i=1}^\infty \frac{|X_i(x)|^2}{\lambda_i}\right)^{\frac{1}{2}} 
\le M_1 \left\|L^2(u_0)\right\| \left( \sum_{i=1}^\infty \frac{|X_i(x)|^2}{\lambda_i}\right)^{\frac{1}{2}}.
$$
Once again,  in the 2nd inequality we used the inequality $\lambda_i^3 < C\lambda_i^4,\ i=1,2,\dots$ with a certain constant $C$. It follows now from the uniform convergence of the 1st of the series in (\ref{4-f6-4}) that  the series (\ref{4-f6-6})  convergences uniformly on any interval $[\epsilon, l],\ 0<\epsilon <l$. The function (\ref{4-f6-1}) is thus $\alpha$-differentiable and the formula
\begin{equation}
\label{4-f6-7}
(D^\alpha_t u)(t) \ = \ - \sum_{i=1}^\infty c_i \lambda_i E_{\alpha}(-\lambda_i t^\alpha)\, X_i(x)
\end{equation}  
holds true. 

We have shown that the function defined by the series (\ref{4-f6})  is at least twice differentiable function with respect to the spatial variable and $\alpha$-differentiable with respect to the time variable and thus is the  solution in the classical sense. \blacksquare

Let us now discuss smoothness of solution (\ref{4-f6}) with respect to the time variable. The following result is valid:

\begin{thm} 
Let $u_0 \in {\cal M}_L$ and $L(u_0) \in {\cal M}_L$. Then the solution (\ref{4-f6}) of the initial-boundary-value problem (\ref{4-ic-h})-(\ref{4-bc-h}) for the one-dimensional time-fractional diffusion equation (\ref{4-eq-h}) is a continuously differentiable function with respect to the time variable on the interval $(0,T)$. 
\end{thm}

Proof:
We already know (see Theorem 5) that (\ref{4-f6}) is the solution of the problem (\ref{4-ic-h})-(\ref{4-bc-h}) for the equation (\ref{4-eq-h}) in the classical sense.  Let us differentiate   the series from the left-hand side of the formula (\ref{4-f6}) with respect to the time variable term by term and construct a series of derivatives:
\begin{equation}
\label{4-f6-8}
\sum_{i=1}^\infty - c_i \lambda_i t^{\alpha -1} E_{\alpha,\alpha}(-\lambda_i t^\alpha)\, X_i(x), \ \ c_i = (u_0, X_i), 
\end{equation}
$E_{\alpha,\beta}$ being the generalized Mittag-Leffler function defined by
\begin{equation}
\label{4-f4}
E_{\alpha,\beta}(z):= \sum_{k=1}^\infty {z^k\over \Gamma(\alpha\,k +\beta)}.
\end{equation}
To get the expression (\ref{4-f6-8}), the formula
$$
\frac{d}{dt} E_{\alpha}(-\lambda_i t^\alpha) = -\lambda_i
t^{\alpha -1} E_{\alpha,\alpha}(-\lambda_i t^\alpha)
$$
was used. The last formula can be obtained by the term by term differentiation of the power series for the Mittag-Leffler function and by using the well-known properties of the Gamma-function. 

Now let $t$ be from the interval $[\epsilon,T],\ 0<\epsilon <T$. Using the estimate (see e.g. \cite{Pod99})
\begin{equation}
\label{MLe-1}
|E_{\alpha,\beta} (-x)|\le {M\over 1+x},\ 0\le x,\ 0<\alpha<1, \ 0<\beta,
\end{equation}
we get the inequality
\begin{equation}
\label{MLe-2}
|\lambda_i t^{\alpha -1} E_{\alpha,\alpha}(-\lambda_i t^\alpha)|\le {M\over t}\frac{\lambda_i t^\alpha}{1+\lambda_it^\alpha}\le \frac{M}{t} \le C,\ \ \epsilon \le t \le T
\end{equation}
with a constant $C$ that of course depends on $\epsilon$ but not on $t$. 

Then we use the estimate (\ref{MLe-2}) and the same reasoning  as in the proof of Theorem 5 and get the following chain of estimates:
$$
\sum_{i=1}^\infty |- c_i \lambda_i t^{\alpha -1} E_{\alpha,\alpha}(-\lambda_i t^\alpha)\, X_i(x)| \le 
C \sum_{i=1}^\infty |c_i\, \lambda_i^{\frac{1}{2}}|\left| \frac{X_i(x)}{\lambda_i^{\frac{1}{2}}}\right| \le 
$$
$$
C_1 \left( \sum_{i=1}^\infty c_i^2\, \lambda_i^2\right)^{\frac{1}{2}} \left( \sum_{i=1}^\infty \frac{|X_i(x)|^2}{\lambda_i}\right)^{\frac{1}{2}} 
\le C_1 \left\|L(u_0)\right\| \left( \sum_{i=1}^\infty \frac{|X_i(x)|^2}{\lambda_i}\right)^{\frac{1}{2}}.
$$
In the 2nd inequality, we used the inequality $\lambda_i < C_2\lambda_i^2,\ i=1,2,\dots$ with a certain constant $C_2$. It follows now from the uniform convergence of the 1st of the series in (\ref{4-f6-4}) that  the series (\ref{4-f6-8})  uniformly converges on any interval $[\epsilon, T],\ 0<\epsilon <T$. The function (\ref{4-f6-1}) is thus continuously differentiable for $0<t<T$ and the formula
\begin{equation}
\label{4-f6-9}
\frac{\partial u}{\partial t} \, =\, -\sum_{i=1}^\infty c_i \lambda_i t^{\alpha -1} E_{\alpha,\alpha}(-\lambda_i t^\alpha)\, X_i(x), \ \ c_i = (u_0, X_i) 
\end{equation}  
holds true. \blacksquare

Now we consider the case of the initial-boundary-value problem (\ref{4-ic-h})-(\ref{4-bc-h}) for the inhomogeneous one-dimensional time-fractional diffusion equation (\ref{4-eq}). 

In the following, we suppose that the initial condition $u_0$ and the function $L(u_0)$ are from the space ${\cal M}_L$. In addition, we suppose that the source function $F$ belongs to the space ${\cal M}_L$ for any $t\in(0,T)$, too.  This means, that for any  $t\in(0,T)$ the function $F$ can be represented in form of an uniformly convergent Fourier series with respect to the spatial variable $x$:
\begin{equation}
\label{4-f6-10}
F(x,t) = \sum_{i=1}^\infty F_i(t)\,X_i(x), \ \ F_i(t) = (F, X_i)=\int_0^l F(x,t)\, X_i(x)\, dx, 
\end{equation} 
$X_i,\ i=1,2,...$ being the eigenfunctions corresponding to the eigenvalues $\lambda_i$
of the eigenvalue problem (\ref{tx4}) - (\ref{tx5}). 

As in the case of the PDEs of parabolic type ($\alpha =1$ in (\ref{4-eq})), we look for solution of the initial-boundary-value problem (\ref{4-ic-h})-(\ref{4-bc-h}) for the equation (\ref{4-eq}) in form of a Fourier series
\begin{equation}
\label{4-f6-11}
u(x,t) = \sum_{i=1}^\infty T_i(t)\,X_i(x), \ \ T_i(t) = (u, X_i).
\end{equation} 
From (\ref{4-f6-11}), the Fourier series for the initial condition (remember the inclusion $u_0\in {\cal M}_L$)
$$
u_0(x) = \sum_{i=1}^\infty c_i\,X_i(x), \ \ c_i = (u_0, X_i),
$$
and the uniqueness of the Fourier series for a function from the space ${\cal M}_L$, we get the initial values for the unknown Fourier coefficients $T_i$:
\begin{equation}
\label{4-f6-12}
T_i(0) = c_i = (u_0, X_i),\ i=1,2,\dots . 
\end{equation} 
Substituting now the representation (\ref{4-f6-11}) into the equation (\ref{4-eq}) and using the properties of the eigenfunctions $X_i$, we get an uncoupled system of initial-value problems for the fractional differential equations for the unknown Fourier coefficients $T_i,\ i=1,2,\dots $:
\begin{equation}
\label{4-f6-13}
\begin{cases}
(D^{\alpha} T_i)(t)=-\lambda_i\, T_i(t) + F_i(t),\\
T_i(0) = c_i.
\end{cases}
\end{equation}

The solution of the initial-value problem (\ref{4-f6-13}) is known (see e.g. \cite{Luc99A}, \cite{Luc99B}) and is given by 
\begin{equation}
\label{4-f6-14}
T_i(t)=c_i\, E_{\alpha}(-\lambda_i t^\alpha)+ \int_0^t \tau^{\alpha -1}
E_{\alpha,\alpha}(-\lambda_i \tau^\alpha)F_i(t-\tau)\, d\tau.
\end{equation}
Thus the formal solution to the initial-boundary-value problem (\ref{4-ic-h})-(\ref{4-bc-h}) for the inhomogeneous one-dimensional time-fractional diffusion equation (\ref{4-eq}) can be written in the form 
\begin{equation}
\label{4-f6-15}
u(x,t)=\sum_{i=1}^\infty \left(c_i\, E_{\alpha}(-\lambda_i t^\alpha)+ \int_0^t \tau^{\alpha -1}
E_{\alpha,\alpha}(-\lambda_i \tau^\alpha)F_i(t-\tau)\, d\tau\right) X_i(x)
\end{equation}
with
$$
c_i = (u_0, X_i),\ \ F_i(t)=(F,X_i).
$$
Using the same technique as in the homogeneous case, it can be shown that 
(\ref{4-f6-15}) is the generalized solution to the problem that can be interpreted as the solution in the classical sense under certain additional conditions. 

Let us now consider asymptotics of the solution (\ref{4-f6-15}) in the case of the homogeneous initial condition $u_0(x)\equiv 0,\ 0\le x \le l$ and for the source function in the form $F(x,t) = \psi(t)\, X_k(x)$, $X_k(x),\ k=1,2,\dots $ being an eigenfunction corresponding to an eigenvalue $\lambda_k$
of the eigenvalue problem (\ref{tx4}) - (\ref{tx5}). The Fourier coefficients of the source function $F$ can be easily evaluated:
$$
F_i(t) =(F,\, X_i)=(\psi(t)X_k(x),\, X_i) = \psi(t)\delta_{ik},
$$
$\delta_{ik}$ being the Kronecker delta symbol. The solution (\ref{4-f6-15}) takes then the form
\begin{equation}
\label{4-f6-16}
u(x,t)= \left( \int_0^t \tau^{\alpha -1}
E_{\alpha,\alpha}(-\lambda_k \tau^\alpha)\psi(t-\tau)\, d\tau\right) X_k(x),
\end{equation}
so we can analyze its behavior for a given function $\psi$. To illustrate the procedure, let us consider the case $\psi(t) = t^{\beta-1}E_{\alpha,\beta}(-\lambda t^{\alpha}),\ \alpha,\beta,\lambda >0$. Then for $\lambda\not = \lambda_k$ the integral in (\ref{4-f6-16}) can be evaluated in a closed form (see \cite{Pod99}) as the Laplace convolution of two Mittag-Leffler functions:
\begin{equation}
\label{4-f6-17}
u(x,t)= \frac{
\lambda E_{\alpha,\alpha+\beta}(-\lambda t^\alpha)-
\lambda_k E_{\alpha,\alpha+\beta}(-\lambda_k t^\alpha)}{\lambda - \lambda_k}
t^{\alpha + \beta -1}X_k(x).
\end{equation}
To extend this result for the case $\lambda = \lambda_k$, we have to differentiate the function $f(\lambda) = \lambda  E_{\alpha,\alpha+\beta}(-\lambda t^\alpha)$ with respect to $\lambda$ and evaluate the result in the point $\lambda = \lambda_k$. The derivative of the Mittag-Leffler function can be determined using the term by term differentiation that leads to the formula ($x\not = 0$)
$$
E^\prime_{\alpha, \beta}(x) = \frac{1}{\alpha x}\left(
E_{\alpha,\beta-1}(x) - (\beta-1)E_{\alpha,\beta}(x)\right).
$$
Using this result, we can explicitly write down the solution $u$ in the case $\lambda = \lambda_k$:
\begin{equation}
\label{4-f6-18}
u(x,t)= \frac{1}{\alpha}\left(
E_{\alpha,\alpha+\beta-1}(-\lambda_k t^\alpha)+ (1-\beta)
E_{\alpha,\alpha+\beta}(-\lambda_k t^\alpha)\right)
t^{\alpha + \beta -1}X_k(x).
\end{equation}
Even if the formulas (\ref{4-f6-17}) and (\ref{4-f6-18}) for the solution $u$ look different in the cases $\lambda \not = \lambda_k$ and $\lambda = \lambda_k$, the asymptotics of solution as $t\to +\infty$ is the same in these two cases. To show this fact, the well-known formula for asymptotics of the Mittag-Leffler function (see e.g. \cite{Pod99})
$$
E_{\alpha,\beta}(-x) = \sum_{k=1}^p \frac{(-x)^{-k}}{\Gamma(\beta-\alpha k)} \, + \, O(x^{-1-p}),\ x\to +\infty,\ p=1,2,\dots .
$$
is used. Applying this formula to the Mittag-Leffler functions in (\ref{4-f6-17}) and (\ref{4-f6-18}), we get the following result:

\begin{thm}
The solution to the initial-boundary-value problem (\ref{4-ic-h})-(\ref{4-bc-h}) for the one-dimensional time-fractional diffusion equation (\ref{4-eq}) with the homogeneous initial and boundary conditions and with the source function in the form 
$$
F(x,t) = \psi(t)\, X_k(x),\ \psi(t)=t^{\beta-1}E_{\alpha,\beta}(-\lambda t^{\alpha}),\ \alpha,\beta,\lambda >0
$$
is given by (\ref{4-f6-17}) in the case $\lambda \not = \lambda_k$ and by (\ref{4-f6-18}) in the case $\lambda = \lambda_k$. The asymptotics of the solution $u(x,t) = \phi(t)X_k(x)$ as $t \to +\infty$ is the same in both cases:
$$
\phi(t)= 
\begin{cases}
O(t^{\beta-\alpha-1}),\ t\to +\infty,\ \mbox{if}\ \alpha \not = \beta, \\ 
O(t^{-\alpha-1}),\ t\to +\infty,\ \mbox{if}\ \alpha = \beta.
\end{cases}
$$
\end{thm}

 
\section{Conclusions and open problems}

In the paper, some initial-boundary-value problems with the Dirichlet boundary conditions for the time-fractional diffusion equation  were considered. Of course, the same method can be applied for the  initial-boundary-value problems with the Neumann, Robin, or mixed boundary conditions. 

The maximum principle enables us to obtain information regarding solutions
of differential equations and  a priori estimates for them without any explicit knowledge
of the form of the solutions themselves, and thus is a valuable tool in scientific
research. In the paper, a maximum principle for the 
time-fractional diffusion equation (\ref{eq}) was discussed and applied for proving  uniqueness of solution of the initial-boundary-value problem  (\ref{ic})-(\ref{bc}) for equation (\ref{eq}). Of course, following the lines of applications of the maximum principle for the parabolic and elliptic PDEs (see e.g. the recent book \cite{MP07}), a lot of other properties of solutions to the time-fractional PDEs can be established. In particular, the maximum principle can be applied for some classes of the non-linear equations of the fractional order, too.

Another important and interesting problem that is still waiting for its solution would be to try to extend the maximum principle to the space- and time-space-fractional PDEs. These equations are very actively employed nowadays in modeling of relevant complex phenomena like anomalous diffusion in inhomogeneous and porous mediums, Levy processes
and Levy flights and the so called fractional kinetics. Similar to the case of time-fractional equations, several different definitions of the space-fractional derivatives are used in these equations. From the physical viewpoint,  an appropriate maximum principle  is expected to be fulfilled for the corresponding models. Clear understanding regarding what definition
of the space-fractional derivatives enables a maximum principle  would essentially help in attempts towards modeling of the real phenomena with  the "right" space- and time-space-fractional PDEs. 

As to the existence of solution to initial-boundary-value problems for the time-fractional diffusion equation, the method presented in the last section can be extended for the n-dimensional problems, too. 



\end{document}